\documentclass[11pt]{article}
\usepackage{amsmath,amssymb,amsfonts,euscript,graphicx}
\usepackage{color}
\usepackage{lmodern}
\usepackage{verbatim}
\usepackage{tikz}
\usetikzlibrary{arrows,shapes}

\usepackage[margin=1.2in]{geometry}

\newtheorem{thm}{Theorem}[section]

\newtheorem{lem}[thm]{Lemma}

\newtheorem{exam}[thm]{Example}

\newtheorem{prop}[thm]{Proposition}

\newtheorem{coro}[thm]{Corollary}

\title{\bf Normal Supercharacter Theory}
\author{Farid Aliniaeifard}

\begin{document}

\maketitle


\begin{abstract}
{  There are two main constructions of supercharacter theories for a group $ G $. The first, defined by Diaconis and 
Isaacs, comes from the action of a group $A$ via automorphisms on our given group $G$.  The second, defined by Hendrickson, is combining a supercharacter theories of a normal subgroup $N$ of $G$ with a supercharacter theory of $G/N$. In this paper we construct a supercharacter theory from an arbitrary set of normal subgroups of $G$. We show that when consider the set of all normal subgroups of G, the corresponding supercharacter theory is related to a partition of $G$ given by certain values on the central idempotents. Also, we show the supercharacter theories that we construct can not be obtained via automorphisms or a single normal subgroup.  }
\end{abstract}


\section{Introduction}

\hspace{4mm} Let $G$ be a finite group, we denote the set of irreducible characters of $G$ by ${\rm Irr}(G)$. The conjugacy class containing $g$ and its cardinality are denoted by $C_g$ and  $m_g$ respectively.  For a subset $S$ of $G$, let $\widehat{S}= \sum_{s\in S}{s}$.\\

Let $U_n(q)$ denote the group of $n \times n$ unipotent upper triangular matrices over a finite field $\mathbb{F}_q$. Classification of the irreducible characters of $U_n(q)$ is a well-known wild problem, provably intractable for arbitrary $n$. In order to find a more tractable way to understand the representation theory of $U_n(q)$, C. Andr\'e \cite{A95} defines and constructs supercharacter theory.  Yan \cite{Y01} shows  how to  replace Andr\'e's construction with more elementary methods.  Diaconis and Isaacs \cite{DI2008} axiomatize the concept of supercharacter theory for an arbitrary group. They mention how to obtain a supercharacter theory for $G$ from the action of $A$ on $G$ by automorphisms.   They also generalize Andr\'e's original construction to define a supercharacter theory for algebra groups, a group of the form $1+J$ where $J$ is a finite dimensional nilpotent associative algebra over a finite field $\mathbb{F}$ of characterestic $p$.  Later, in \cite{H12}, Hendrickson shows how to construct other supercharacter theories for an arbitrary group $G$ by combining certain supercharacter theory for a normal subgroup $N$ of $G$ with a supercharacter theory for $G/N$. 
Also in \cite{N12} the authors  obtain a relationship
between the supercharacter theory of all unipotent upper triangular matrices over a finite field $\mathbb{F}_q$ simultaneously and
the combinatorial Hopf algebra of symmetric functions in non-commuting variables.
\\

Let $N(G)$ be the set of all normal subgroup of $G$. Let $\{N_1,\ldots,N_k\}\subseteq N(G)$.  We define
$A(N_1,\ldots,N_k)$ to be the smallest subset of $N(G)$  such that\\
\\
(1) $e,G\in A(N_1,\ldots,N_k)$.\\
(2) $\{N_1,\ldots,N_k\}\subseteq A(N_1,\ldots,N_k)$.\\
(3) $A(N_1,\ldots,N_k)$ is closed under product and intersection. \\
\\
Define $${N^{\circ}_{A(N_1,\cdots,N_k)}}=N\setminus \bigcup_{K\in A(N_1,\cdots,N_k), K< N} K.$$ 
For simplicity of notation, we write ${N^{\circ}}$ instead of ${N^{\circ}_{A(N_1,\ldots,N_k)}}$ when it is clear that $N$ is in $A(N_1,\ldots,N_k)$. We will show that $\{N^{\circ}: N\in A(N_1,\ldots,N_k)\}$ is the set of superclasses of a supercharacter theory,  and we call such supercharacter theory the {\it normal superchracter theory generated by $\{N_1,\ldots,N_k\}$}.  In general this supercharacter theories can not be constructed by the previous supercharcter theory constructions. Remark that when we have a larger set of normal subgroups, the normal supercharacter theory we obtain will be finer. In particular the finest normal supercharacter theory is obtained when we consider the set of all normal subgroups of $G$, and is related to a partition of $G$ given by certain values on the central idempotents.\\

In Section 2, we review definitions and notations for supercharacter theories. In particular we mention the known constructions of supercharacter theories.
Next in Section 3, we define our normal supercharacter theory construction. We also show that the finest normal supercharacter theory is obtained by considering certain values of the central idempotents. In Section 4,  we show that the normal supercharacter theory  can not be obtained by the previous general constructions. Finally, in the last section we list some open problems related to the concept.

\section*{Acknowledgment}
\hspace{4mm}  I wish to express my appreciation to my supervisor professor Bergeron who
carefully read an earlier version of this paper and made significant suggestions for
improvement. Also, I would like to thank Shu Xiao Li for his helpful comments.

\section{Background } 

 \hspace{4mm} We first mention the definition of supercharacter theory by Diaconis and Isaacs \cite{DI2008}.\\

A supercharacter theory of a finite group $G$ is a pair ($\mathcal{X}$,$\mathcal{K}$) where 
$\mathcal{X}$ is a partition of  ${\rm Irr}(G)$  and $\mathcal{K}$ is a partition of $G$ such that:
\\
\\
(a)  $|\mathcal{K}|=|\mathcal{X}|$,\\
(b)   for $X \in \mathcal{X} $, the character $\mathcal{X}_X$, a nonzero character whose irreducible constituents lie in $X$, is constant on the parts of $\mathcal{K}$ and\\
(c)  the set $\{1\}\in \mathcal{K}$.\\
\\
We will refer to characters $\mathcal{X}_X$ as the supercharacters and  to the member of $\mathcal{K}$ as superclasses.\\

Every finite group has two trivial supercharacter theories: the usual irreducible character theory
and the supercharacter theory $(\{{\bf 1}\},\{{\rm Irr}(G) \setminus {\bf 1}\}\},\{ \{1\}, G \setminus \{1\}\} )$, where
{\bf 1} is the principal character of $G$. \\

The concept of a Schur ring is defined by Schur in \cite{S1933}.  Hendrickson \cite{H12} shows that there is a bijection between the supercharacter theories of a group $G$ and Schur rings over $G$ contained in $Z(\mathbb{C}[G])$, the center of $\mathbb{C}[G]$.\\
\\
{\bf Definition.} Let $G$ be a finite group. A subring $S$ of the group algebra $\mathbb{C}[G]$ is called a Schur ring over $G$ if there is a set partition $\mathcal{K}$  of $G$ such that $\{1\}\in \mathcal{K}$, $S=\mathbb{C}$-span$\{\widehat{K}: K\in \mathcal{K}\}$, and $\{g^{-1}: g\in K\}\in \mathcal{K}$ for all $K\in \mathcal{K}$.

\begin{prop}{\rm \cite[Proposition 2.4]{H12}}
Let $G$ be a finite group. Then there is a bijection
$$\left\lbrace Supercharacter~ theories~(\mathcal{X},\mathcal{K})~of~G  \right\rbrace \longleftrightarrow \left\lbrace  Schur~ rings~ over~G~ contained ~in~ Z(\mathbb{C}[G])  \right\rbrace \newline$$
$~~~~~~~~~~~~~~~~~~~~~~~~~~~~~(\mathcal{X},\mathcal{K})~~~~~~~~~~~~~\longmapsto~~~~~~~~~~~ \mathbb{C}$-{\rm span}$\{\widehat{K}:K\in \mathcal{K}\}.$

\end{prop}

In the proof of surjectivity of the above bijection, Hendrickson does not need the condition $\{g^{-1}: g\in K\}\in \mathcal{K}$ from the definition of Schur ring. So we have the following corollary.

\begin{coro}\label{sup}
Let $G$ be a finite group and let  $\mathcal{K}$ be a partition of $G$. Then the following statements are equivalent.

{\rm (1)} $\mathcal{K}$ is the set of superclasses of a supercharacter theory.

{\rm (2)}  $\{1\}\in \mathcal{K}$ and $\mathbb{C}$-{\rm span}$\{\widehat{K}:K\in \mathcal{K}\}$ is a subring of $Z(\mathbb{C}[G])$.
\end{coro}
{\bf Definition.} A superclass theory is a partition $\mathcal{K}$ of $G$ satisfying one of the two equivalent conditions in Corollary \ref{sup}.\\

Now we discuss two main methods of constructing 
supercharacter theories of an arbitrary finite group.\\

\subsection{A Group Acts Via Automorphisms on a Given Group}

\hspace{4mm} Given finite groups $A$ and $G$, we say that $A$ acts via automorphisms  on $G$ if $A$ acts on $G$ as a set, and in addition $(gh).x=(g.x)(h.x)$ for all $g,h\in G$ and $x\in A$. An action via automorphisms of $A$ on $G$ determines and is determined by a homomorphism $\phi: A\rightarrow Aut(G)$.\\

 Suppose that $A$ is a group that acts via automorphisms on our given group $G$. It is well known that
$A$ permutes both the irreducible characters of $G$ and the conjugacy classes of $G$. By a
lemma of R. Brauer, the permutation characters of $A$ corresponding to these two actions
are identical, and so the numbers of $A$-orbits on ${\rm Irr}(G)$ and on the set of
classes of $G$ are equal (See Theorem 6.32 and Corollary 6.33 of \cite{I1994}).  It is easy to see that these orbit
decompositions yield a supercharacter theory $(\mathcal{X},\mathcal{K})$ where members of $\mathcal{X}$ are the  $A$-orbits on ${\rm Irr}(G)$ and members of $\mathcal{K}$ are the unions of the $A$-orbits on the classes of $G$. It is
clear that in this situation, the sum of the characters in an orbit $X \in \mathcal{X}$ is constant on
each member of $\mathcal{K}$. We denote by ${\rm AutSup}(G)$ the set of all such supercharacter theories of $G$.

\subsection{$\ast$-Product}

\hspace{4mm} Suppose that $A$ is a group that acts via automorphisms on our given group $G$. Let  ${\rm Sup}(G)$ be the set of all supercharacter theories of $G$. We say that $(\mathcal{X},\mathcal{K})\in {\rm Sup}(G)$ is $A$-invariant if the action of $A$ fixes each part $K \in\mathcal{K}$ setwise. We denote by ${\rm Sup}_{A}(G)$ the set of $A$-invariant supercharacter theories of $G$.
Note that  if $N$ is normal in $G$, then $\mathcal{C} \in {\rm Sup}(H)$ is $G$-invariant if and only if its superclasses are union of conjugacy classes of $G$. Also, if $M,N$ are normal subgroup of $G$ and $N<M$, then a supercharacter theory of $M/N$ is $G/N$-invariant if and only if it is $G$-invariant.\\
\\
{\bf Notation.} Let $N$ be a normal subgroup of a group $G$. If $\mathcal{L}$ is a set of subsets of $G/N$, then we define 
$\widetilde{\mathcal{L}}=\{\cap_{Ng\in L}Ng : L\in \mathcal{L}\}$. Let $\chi\in {\rm Irr}(N)$. We denote  by  ${\rm Irr}(G|\psi)$ the set of irreducible characters $\psi$ of $G$ such that the inner product of $\psi$ and $\chi$ is positive.  If $\mathcal{Z}$ is a set of subsets of ${\rm Irr}(N)$, then we define $\mathcal{Z}^{G}=\{\cup_{\psi \in Z}{\rm Irr}(G|\psi):Z\in \mathcal{Z}\}$.  Now consider $(\mathcal{X},\mathcal{K})\in {\rm Sup}_G(N)$. Since $\{1_N\}\in \mathcal{X}$, one part of $\mathcal{X}^{G}$ is $\{1_N\}^{G}=\{\chi\in {\rm Irr}(G): N\subseteq {\rm ker} \chi\}$, which we identify with ${\rm Irr}(G/N)$ in the usual natural way.

\begin{thm}{\rm \cite[Theorem 4.3]{H12}}
 Let $G$ be a group and $N$ be a normal subgroup of $G$. Let $\mathcal{C}=(\mathcal{X},\mathcal{K})\in {\rm Sup}_G(N)$ and $\mathcal{D}=(\mathcal{Y},\mathcal{L})\in {\rm Sup}(G/N)$. Then  
$$(\mathcal{Y}\cup \mathcal{X}^{G}\setminus\{{\rm Irr}(G/N)\}, \mathcal{K}\cup \widetilde{\mathcal{L}}\setminus \{N\})$$ is a supercharacter theory of $G$.
\end{thm}

We call the supercharacter theory of $G$ constructed in the procending theorem the $\ast$-product of $(\mathcal{X},\mathcal{K})$ and $(\mathcal{Y},\mathcal{L})$,  and write it as $(\mathcal{X},\mathcal{K}) \ast (\mathcal{Y},\mathcal{L})$. Also, let ${\rm Sup}^{*}(G)$  denote the set of all supercharacter theories of $G$ obtained by $*$-product.

\section{Normal Supercharacter Theory}

\hspace{4mm} In this chapter we construct a supercharacter theory from an arbitrary set of normal subgroups. We call such supercharacter theory  a normal supercharacter theory.  

\subsection{Supercharacter Theory From Central Idempotents}

\hspace{4mm} In this section, we consider a partition of conjugacy classes and irreducible characters given by certain values of central idempotent. In the next section we will see that it is a supercharacter theory and is given by the finest normal supercharacter theory.\\

 By \cite[Proposition 8.15]{L1991} every character $\chi \in {\rm Irr}(G)$ has a corresponding central idempotent $$e_\chi=|G|^{-1} \chi(1) \sum_{g\in G}{\chi(g^{-1})g}.$$ These idempotents are orthogonal, i.e, $e_\chi e_\phi=0$ when $\chi\neq \phi$. Recall that $\widehat{C_g}=\sum_{i}{m_g {\chi_i (1)}^{-1}\chi_{i}(g) e_{\chi_i}}$. Therefore, 

$$m_g1-\widehat{C_g}=m_g1-\sum_{i}{m_g {\chi_i (1)}^{-1}\chi_{i}(g) e_{\chi_i }}=m_g(1-\sum_{i}{ {\chi_i (1)}^{-1}\chi_{i}(g) e_{\chi_i}})=$$ 

$$m_g(\sum_{i}{e_{\chi_{i}}}-\sum_{i}{ {\chi_i (1)}^{-1}\chi_{i}(g) e_{\chi_i }})=m_g(\sum_{i}{ (1-{\chi_i (1)}^{-1}\chi_{i}(g)) e_{\chi_i }})$$

$$\Rightarrow 1-\frac{\widehat{C_g}}{m_g}=\sum_{i}{ (1-\frac{\chi_{i}(g)}{{\chi_i (1)}}) e_{\chi_i }}=\sum_{i}{ (1-\frac{\chi_{i}(g)}{{\chi_{i}(1)}}) e_{\chi_i}}.$$
Look at the last equation i.e., $1-\frac{\widehat{C_g}}{m_g}=\sum_{i}{ (1-\frac{\chi_{i}(g)}{{\chi_{i}(1)}}) e_i}$. Let 

$$E_g=\{e_{\chi_i}: 1-\frac{\chi_{i}(g)}{{\chi_{i}(1)}}\neq 0\}, K_g=\cup_{E_g=E_h}{C_h}, ~{\rm and}~ U_g=\cup_{E_h\subseteq E_g}{C_h}.$$ 


As in following example we will see that $\{K_g : g\in G\}$ is a superclass theory.

\begin{exam}
The character table of $S_5$ is

{\rm
$$\begin{tabular}{l*{6}{c}r}
Classes             & (1) & (1~2) & (1~2~3) & (1~2~3~4) & (1~2~3~4~5)  & (1~2)(3~4) & (1~2)(3~4~5) \\
\hline
$\chi_1$            &1 & 1 & 1 & 1 & 1 & 1 & 1~~~~~~ \\
$\chi_2$            & 1 & -1 & 1 & -1 &  1 & 1 &-1~~~~~~ \\
$\chi_3$             & 4& 2 & 1 & 0 &  -1 & 0& -1 ~~~~~~\\
$\chi_4$           & 4 & -2 & 1 & 0 &  -1 & 0 &  1~~~~~~ \\
$\chi_5$           & 5 & -1 & -1& 1 &  0& 1 &  -1~~~~~~ \\
$\chi_6$           & 5 & 1& -1& -1&  0& 1&  1~~~~~~ \\
$\chi_7$           & 6 & 0 & 0 & 0 &  1 & -2 &  0~~~~~~  
\end{tabular}$$
}
by definition 

$$E_{(1~2~3)}=E_{(1~2~3~4~5)}=E_{(1~2)(3~4)},$$

 then 

$$K_{(1~2~3)}=K_{(1~2~3~4~5)}=K_{(1~2)(3~4)}.$$

Also, 

$$E_{(1~2)}=E_{(1~2~3~4)}=E_{(1~2)(3~4~5)},$$ 

then $$K_{(1~2)}=K_{(1~2~3~4)}=K_{(1~2)(3~4~5)}.$$

One can check that $(\{\{e\},K_{(1~2~3)},K_{(1~2)}\},\{\{\chi_1\},\{\chi_2\}, \{\chi_3,\chi_4,\chi_5,\chi_6,\chi_7\} \})$ is a supercharacter theory. \hfill $\square$

\end{exam}

In the above example $\{K_g:g\in G\}$ forms a superclass theory. A natural question arises: does $\{K_g:g\in G\}$ always give rise to a superclass theory? We will answer this question in Corollary \ref{kg}.

\subsection{Normal Supercharacter Theory}

\hspace{4mm} In this section we construct our normal supercharacter theory. We will show the finest normal supercharacter theory is related to $\{K_g:g\in G\}$ the partition of $G$ given by $\{E_g: g\in G\}$ a subset of the set of all subsets of central idempotents. We need the following definitions and notations in the sequel.\\

If $(P, \leq)$ is a poset and $\Bbb{C}^{P\times P}$ is the set of all functions $\alpha: P\times P \rightarrow \Bbb{C}$, the associated incidence algebra is
$$A(P) = \{\alpha \in \Bbb{C}^{
P \times P}
: \alpha(x, y) = 0 ~{\rm unless}~ x \leq y\}.
$$
The mobius function $\mu\in A(P)$ is defined recursively  by the following rule: $\mu(x, y) = 0$ whenever
$x \not \leq y$, $\mu(x, x) = 1$ and for $x<y$
$$\mu(x, y) =- \sum_{x< z\leq y}{\mu(z,y)}.$$
It is immediate from this definition that
\[  \sum_{x\leq z\leq y}{\mu(z,y)}= \left\{ 
  \begin{array}{l l}
    1 & \quad \text{if $x=y$ }\\
    0 & \quad \text{otherwise}
  \end{array}
	\right.\]

Let $N(G)$ be the set of all normal subgroup of $G$. Note that the product of two normal subgroup is a normal subgroup. We can see that  $N(G)$ is a semigroup. Let $\{N_1,\ldots,N_k\}\subseteq N(G)$.  We define
$A(N_1,\ldots,N_k)$ to be the smallest subsemigroup of $N(G)$ containing $\{N_1,\ldots,N_k\}$ such that\\
\\
(1) $e,G\in A(N_1,\ldots,N_k)$.\\
(2) $A(N_1,\ldots,N_k)$ is closed under  intersection. \\

 Note that every element $N\in A(N_1,\ldots,N_k)$ is a normal subgroup of $G$. We define for an element $N\in A(N_1,\ldots,N_k)$
 $${N^{\circ}_{A(N_1,\cdots,N_k)}}=N\setminus \bigcup_{K\in A(N_1,\cdots,N_k), K< N} K.$$ 
For simplicity of notation, we write ${N^{\circ}}$ instead of ${N^{\circ}_{A(N_1,\ldots,N_k)}}$ when it is clear that $N$ is in $A(N_1,\ldots,N_k)$.
 Note that    $\widehat{{N}}=\sum_{ H\in A(N_1,\ldots,N_k),\{e\}\leq H\leq N }\widehat{H^{\circ}}$. Thus, by Mobius Inversion Theorem we have 
 $$\widehat{{N^{\circ}}}=\sum_{ H\in A(N_1,\cdots,N_k),\{e\}\leq H\leq N }{\mu(H,N)}\widehat{H}.$$

\begin{exam}
Let $G=C_2\times C_4$.  Here is the Hasse diagram for $N(G)$, and $\mu(H,G)$ for every $H\in N(G)$ is written above the vertex $H$.

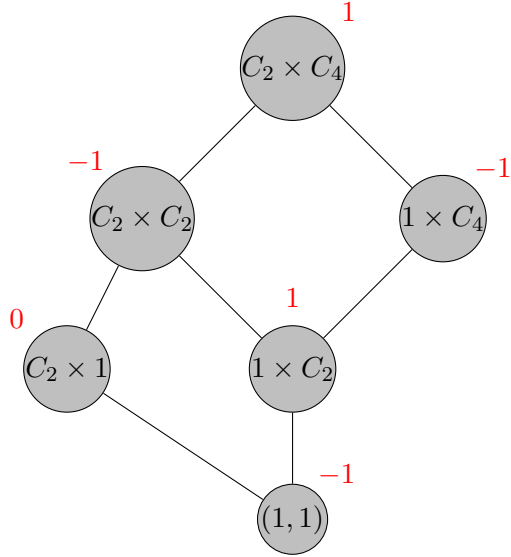
\begin{center}
\begin{tikzpicture}[scale=1, auto,swap]
\tikzstyle{every node}=[circle,fill=black!25,minimum size=20pt,inner sep=0pt];
\node (a) at (0,0)[draw, circle, fill][label=above right :$\color{red}{-1}$]{$ (1,1)$};
\node (b) at (-3,2)[draw, circle, fill][label=above left:$\color{red}{0}$]{$ C_2\times 1$};
\node (c) at (0,2)[draw, circle, fill][label=above:$\color{red}{1}$]{$ 1 \times C_2$};
\node (d) at (-2,4)[draw, circle, fill][label=above left:$\color{red}{-1}$]{$ C_2\times C_2$};
\node (e) at (2,4)[draw, circle, fill][label=above right:$\color{red}{-1}$]{$1\times C_4$};
\node (f) at (0,6)[draw, circle, fill][label=above right:$\color{red}{1}$]{$C_2\times C_4$};

\draw (a)--(b)(a)--(c)--(d)--(b)(e)--(c)(d)--(f)--(e);
																		
\end{tikzpicture}

\end{center} 
By the above diagram it is easy to see that $\widehat{{C_2\times C_4}^{\circ}}=\sum_{\{e\}\leq H\leq C_2\times C_4}{\mu(H,C_2\times C_4)}\widehat{H}.$         $\hfill \square$

\end{exam}

\begin{thm}\label{normal}
Let $\{N_1,\ldots,N_k\}$ be  a set of normal subgroups of $G$.  Then $\mathcal{K}=\{{N^{\circ}}: N\in A(N_1,\ldots,N_k)\}$ is a superclass theory. 
\end{thm}
{\bf Proof.} By the definition it is clear that if $N^{\circ},H^{\circ} \in \mathcal{K}$, then $N^{\circ}=H^{\circ}$ or  $N^{\circ}\cap H^{\circ}=\emptyset$. Furthermore, $\bigcup_{N\in A(N_1,\ldots,N_k)}{N^{\circ}}=G$. Therefore, $\mathcal{K}$ is a partition of $G$ such that $\{e\}\in \mathcal{K}$.

 Let $g\in {N^{\circ}}$. Then $C_g\in N$. If $C_g\cap H\neq \emptyset$ for some normal subgroup $H$ in $A(N_1,\ldots,N_k)$ and $H\subset N$, then $g\in C_g\subseteq H$. Thus, $g\not \in {N^{\circ}}$,  yielding a contradiction. Therefore, we must have $C_g\subseteq {N^{\circ}}$. So every member of $\mathcal{K}$ is a union of conjugacy classes of $G$. 
 We have
 $$\widehat{{N^{\circ}}}=\sum_{ H\in A(N_1,\cdots,N_k), H\leq N}{\mu(H,N)}\widehat{H}~~~{\rm and} ~~~\widehat{N}=\sum_{H\in A(N_1,\cdots,N_k), H\leq N} {\widehat{H^{\circ}}}.$$
 Therefore, $\mathbb{C}$-span$\{\widehat{N^{\circ}}: N\in A(N_1,\ldots,N_k)\}=\mathbb{C}$-span$\{\widehat{N}: N\in A(N_1,\ldots,N_k)\}$. Since $A(N_1, \ldots, N_k)$  is closed under product, $\mathbb{C}$-span$\{\widehat{N}: N\in A(N_1,\ldots,N_k)\}$ is a subalgebra of $Z(\mathbb{C}[G])$. We conclude by Corollary \ref{sup} that $\mathcal{K}$ is a superclass theory.
 \hfill $\square$
\\
\\
As you see in Theorem \ref{normal}, for a set of normal subgroups $\{N_1,\ldots,N_k\}$ of $G$, $\{{N^{\circ}}: N\in A(N_1,\ldots,N_k)\}$ is a superclass theory. We say a supercharacter theory $(\mathcal{X},\mathcal{K})$ is a normal supercharacter theory if $\mathcal{K}=\{{N^{\circ}}: N\in A(N_1,\ldots,N_k)\}$ for some normal subgroups $N_1,\cdots, N_k$ of $G$. We denote by $NSup(G)$ the set of all possible normal supercharacter theories of $G$. 
\\
\\

A subgroup of $G$ is normal if and only if it is the union of a set of conjugacy classes of $G$. We have
an equivalent characterization of normality in terms of the kernels of irreducible characters. Recall
that the kernel of a character $\chi$ of $G$ is the set ${\rm ker} \chi = \{g \in G: \chi(g) = \chi(1)\}$. This is just the kernel
of any representation whose character is $\chi$, and so ${\rm ker} \chi$ is normal subgroup. A subgroup of $G$ is
normal if and only if it is the intersection of the kernels of some finite set of irreducible characters
\cite[Proposition 17.5]{JL1993}; thus the normal subgroups of $G$ are the subgroups which we can construct
from the character table of $G$. 

Recall that 
$$E_g=\{e_{\chi_i}: 1-\frac{\chi_{i}(g)}{{\chi_{i}(1)}}\neq 0\}~{\rm and} ~K_g=\cup_{E_g=E_h}{C_h}.$$ We show that $K_g=N^{\circ}$ for a normal subgroup of $G$, and if for a normal subgroup $N$ of $G$, $N^{\circ}\neq \emptyset$, then there is a $g\in N$ such that $K_g=N^{\circ}$. First we prove the following lemma.

\begin{lem}\label{kgg}
Assume that $E_g=\{e_{\chi_{i_1}},\cdots,e_{\chi_{i_k}}\}$. If 
$$N= \bigcap_{\chi \in {\rm Irr}(G) \setminus \{\chi_{i_1},\cdots,\chi_{i_k}\}}{\rm ker}(\chi),$$
then $K_g= {N^{\circ}}$.
\end{lem}
{\bf Proof.} 
Let $h\in K_g$. Then $
E_g=E_h$, and  so for every $\chi  \in {\rm Irr}(G)\setminus  \{\chi_{i_1},\cdots,\chi_{i_k}\}$, $\chi(h)=\chi(1)$. Therefore,
 $$h\in \bigcap_{\chi \in {\rm Irr}(G)\setminus \{\chi_{i_1},\cdots,\chi_{i_k}\}}{\rm ker}(\chi)=N.$$
 Let $H$ be a normal 
subgroup of $G$ such that $H\subset N$. Since every normal subgroup is a intersection of some kernels of irreducible characters and  $H\subset N$, 
we have $$H=(\cap_{\chi \in {\rm Irr}(G)\setminus\{\chi_{i_1},\cdots,\chi_{i_k}\}}{\rm ker}(\chi) )\bigcap (\chi_{j_1}\cap \cdots \cap \chi_{j_t})$$ where
 $\chi_{j_1}, \cdots, \chi_{j_t}\in \{\chi_{i_1},\cdots,\chi_{i_k}\}\neq \emptyset$. If $h\in H$, then $h\in {\rm ker}(\chi_{j_1})\cap  \cdots \cap {\rm ker}(\chi_{j_t})$, i.e., $h\not \in K_g$, yielding a contradiction.
  Therefore, $K_g\cap H=\emptyset$ for every normal subgroup $H\subset N$. Thus, $h\in N^{\circ}$ and so $K_g\subseteq N^{\circ}$. Now let $h\in N^{\circ}$. Then $E_g\subseteq E_h$. If $E_g\neq E_h$, there is a irreducible character $\chi\in  \{\chi_{i_1},\cdots,\chi_{i_k}\} $ 
  such that $\chi(h)=\chi(1)$. Then $g\in H=\bigcap_{\psi\in {\rm Irr}(G)\setminus  \{\chi_{i_1},\cdots,\chi_{i_k}\}}{\rm ker}(\psi) \cap {\rm ker}(\chi)$. Thus, $h\in H\subset N$, yielding a contradiction. We conclude that $E_g=E_h$, i.e., $h\in K_g$, and so $N^{\circ} \subseteq K_g.$
	$\hfill \square$

\begin{thm}
Let $G$ be a group. Then
\\
\\
{\rm(1)} For every $g\in G$, $K_g={N^{\circ}}$ for some normal subgroup $N$ of $G$.\\
{\rm (2)} For every normal subgroup $N$ of $G$ such that ${N^{\circ}}\neq \emptyset$, there is an element $g\in {N^{\circ}}$ such that $K_g={N^{\circ}}$.

\end{thm}
{\bf Proof.}
(1) Let $E_g=\{e_{\chi_{i_1}},\cdots,e_{\chi_{i_k}}\}$ and let 
$$N= \bigcap_{\chi \in {\rm Irr}(G) \setminus \{\chi_{i_1},\cdots,\chi_{i_k}\}}{\rm ker}(\chi).$$
Then by Lemma \ref{kgg}, $K_g= {N^{\circ}}$. 
\\
\\
(2) Let $N$ be a normal subgroup of $G$ such that ${N^{\circ}}\neq \emptyset$. Let $g\in {N^{\circ}}$. We show that $K_g={N^{\circ}}$. Since $N$ is a normal subgroup, $N=\bigcap_{i\in I}{\rm ker}(\chi_{i})$.  If there is an irreducible character $\chi\in {\rm Irr}(G)\setminus \{\chi_i:i\in I\}$ such that $\chi(g)=\chi(1)$, then $g\in H=\bigcap_{i\in I}{\rm ker}(\chi_{i}) \cap {\rm ker}(\chi)$. Thus, $g\in H\subset N$, yielding a contradiction. Therefore, $E_g=\{e_{\chi} : \chi \in {\rm Irr}(G) \setminus\{\chi_i: i \in I\}\}$.
 So by Lemma \ref{kgg}, $K_g={N^{\circ}}$.  
\hfill $\square$

\begin{coro}
Let $G$ be a group. Then  for every $g\in G$, $U_g$ is a normal subgroup of $G$.
\end{coro}
{\bf Proof.} Recall that $$ U_g=\cup_{E_h\subseteq E_g}{C_h}.$$
Let $E_g=\{e_{\chi_{i_1}},\cdots,e_{\chi_{i_k}}\}$. We show that $$U_g=N= \bigcap_{\chi \in {\rm Irr}(G)\setminus \{\chi_{i_1},\cdots,\chi_{i_k}\}}{\rm ker}(\chi).$$
Note that $N=\bigcup_{K\in N(G),K\leq N}{{K^{\circ}}}$. Let $K\in N(G)$ and  $K\subset N$. If ${K^{\circ}}=\emptyset$, then  ${K^{\circ}}\subseteq U_g$. If
 ${K^{\circ}}\neq \emptyset$, then by Proposition \ref{kg}, there is a $h\in K^{\circ}$ such that $K_h={K^{\circ}}$. Since $h\in K\subseteq N$, we have
  $\chi(h)=\chi(1)$ for every   ${\chi \in {\rm Irr}(G)\setminus \{\chi_{i_1},\cdots,\chi_{i_k}\}}$. Therefore, $E_h\subseteq E_g$, and so
  $K^{\circ}=K_g\subseteq U_g$. We conclude that $N\subseteq U_g$.  Now we want to show that $U_g\subseteq N$. Let $h\in U_g$. Then $E_h\subseteq E_g$, and
 so $h\in \bigcap_{\chi \in {\rm Irr}(G) \setminus \{\chi_{i_1},\cdots,\chi_{i_k}\}}{\rm ker}(\chi)=N$. Therefore, $U_g\subseteq N$.
\hfill $\square$\\

As we mentioned before, the finset normal supercharacter theory is when we generate a normal supercharcter theory by $N(G)$ the set of all normal subgroups of $G$. In the following corollary we show that the finest normal supercharacter theory is equal to the supercharacter theory with $\{K_g : g\in G\}$ as the set of superclasses. And since every $K_g$ is related to a set of central idempotents, we can see that the finest normal supercharacter theory corresponds to  a set of subsets of central idempotents.

\begin{coro}\label{kg}
Let $G$ be a group. Then $\{K_g : g\in G\}=\{{N^{\circ}} : N\in N(G)\}$ is the finest superclass theory.
\end{coro}
{\bf Proof.}
The normal supercharacter theory generated by $N(G)$ has $\{{N^{\circ}} : N\in N(G)\}$ as the set of superclasses. Since every non-empty ${N^{\circ}}$ is equal to $K_g$ for some $g\in G$. We have $\{{N^{\circ}} : N\in N(G)\}=\{K_g:g\in G\}$. Therefore, $\{K_g : g\in G\}$ is a superclass theory. \hfill $\square$



\section{$NSup(G)$ is not a subset of the union of $Sup^{*}(G)$ and $AutSup(G)$}

\hspace{4mm} In the following example we show that $Sup^{*}(G)\cap AutSup(G)$ is not a subset of $NSup(G)$ and there is a normal supercharacter theory which is not in the union of $Sup^{*}(G)$ and $AutSup(G)$.

\begin{exam} Let $G=C_3 \times C_4$. Note that the supercharacter theory correspond to superclass theory $\{ C_g: g\in G\}$  is in $Sup^{*}(G)\cap AutSup(G)$, but it is not in $NSup(G)$. Therefore, $Sup^{*}(G)\cap AutSup(G)$ is not a subset of $NSup(G)$.   We now construct normal supercharacter theory generated by $\{C_3 \times 1, 1\times C_4\}$ and we show that it is not  in the union of $Sup^{*}(G)$ and $AutSup(G)$.

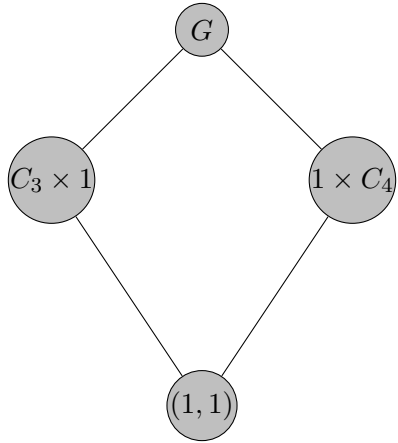
\begin{center}
\begin{tikzpicture}[scale=1, auto,swap]
\tikzstyle{every node}=[circle,fill=black!25,minimum size=20pt,inner sep=0pt];
\node (a) at (0,1)[draw, circle, fill][label=right:${}$]{$ (1,1)$};
\node (b) at (-2,4)[draw, circle, fill][label=left:${}$]{$ C_3\times 1$};
\node (c) at (2,4)[draw, circle, fill][label=right:${}$]{$1\times C_4$};
\node (d) at (0,6)[draw, circle, fill][label=right:${}$]{$G$};

\draw (a)--(b)(c)--(d)(b)--(d)(a)--(c);
																		
\end{tikzpicture}

\end{center}
The set of superclasses for normal supercharacter theory generated by  $\{C_3 \times 1, 1\times C_4\}$ is $\{1, \{(g,1),(g^2,1)\},\{(1,h),(1,h^2),(1,h^3)\},\{(g,h),(g,h^2), (g,h^3),(g^2,h),(g^2,h^2),(g^2,h^3)\}\}$.\\

 Let $(\mathcal{X},\mathcal{K})\in AutSup(G)$. Since $Aut(G)\cong \Bbb{Z}_2 \times \Bbb{Z}_2$, every $A$-orbit has at most $4$ members. Note that the members of $\mathcal{K}$ are the unions of the $A$-orbits on the classes of $G$. Therefore, every members of $\mathcal{K}$ has at most cardinality $4$. But we have a superclass with cardinality $6$ in the normal  supercharacter theory generated by $\{C_3 \times 1, 1\times C_4\}$. Thus, this normal supercharacter theory for $G$ is not in $AutSup(G)$.
\\

 Now we show that normal supercharacter theory generated by  $\{C_3 \times 1, 1\times C_4\}$ is not in $Sup^{*}(G)$. If we choose a subgroup of order $2$ and construct the $\ast$-product, then there is two superclass with cardinality 1, but we just have one superclass with cardinality $1$ in the normal supercharacter theory generated by  $\{C_3 \times 1, 1\times C_4\}$. Let us choose a subgroup of order $4$. Then $\{(g,1),(g^2,1)\}$ is not a superclass of this supercharacter theory. 
Now we choose a subgroup of order 3, and construct the supercharacter theory by $\ast$-product. Then $\{(1,h),(1,h^2),(1,h^3)\}$ is not a superclass of this supercharacter theory.
\\

Therefore, the normal supercharacter theory generated by  $\{C_3 \times 1, 1\times C_4\}$ is not in the union of $Sup^{*}(G)$ and $AutSup(G)$.

\end{exam}

\section{Related Problems}

\hspace{4mm} There are the following open problems related to this subject.\\
\\
(1) A supercharacter theory  $(\mathcal{X}, \mathcal{K})$ is said to be integral if $X(g)$ is an integer for every $X\in \mathcal{X}$ and $g\in G$. For which group $G$ and which set of normal subgroups of $G$ the normal supercharacter theory is integral.
\\
\\
(2) In \cite{N12}, there is a Hopf algebra structure for all supercharacter theories of unipotent uppertriangular matrices over a finite field which is isomorphic to Hopf algebra structure of symmetric functions on noncommutative variables. Is there a Hopf  Algerba structure for the set of a normal supercharacter theories of all unipotent uppertriangular matrices over a finite field? \\
\\
(3) Does the supercharacter theory for group algebras defined in \cite{DI2008} a normal supercharacter theory?\\
\\
(3) What is the relation between normal supercharacter theory of a group algebra and supernormal subgroups with respect to that supercharacter theory.
\\
\\
(5) What are the supercharacters for a normal supercharacter theory?
\\
\\
(6) Can we unify the available supercharacter theories?

Farid Aliniaeifard, Department of Mathematics and Statistics, York University, Toronto, Ontario, Canada M3J 1P3.
 E-mail address: faridand@mathstat.yorku.ca


\begin{thebibliography}{1}
 
 \bibitem{N12}M. Aguiar, C. Andr\'e, C. Benedetti, N. Bergeron, Z. Chen, P. Diaconis, A. Hendrickson, S.
 Hsiao, M. Isaacs, A. Jedwab, K. Johnson, G. Karaali, A. Lauve, T. Le, S. Lewis, H. Li, K.
 Magaard, E. Marberg, J-C. Novelli, A. Pang, F. Saliola, L. Tevlin, J-Y. Thibon, N. Thiem, V.
 Venkateswaran, C.Vinroot, N. Yan and M. Zabrocki, Supercharacters, symmetric functions in
 noncommuting variables, and related Hopf algebras. Adv. Math. 229, Issue 4. (2012).\vspace{-3mm}
 \bibitem{A95} C. Andr\'e, Basic characters of the unitriangular group, J. Algebra 175 (1995) 287-319.\vspace{-3mm}
 \bibitem{A98} C. Andr\'e, Irreducible characters of finite algebra groups, in: Matrices and Group Representations, Coimbra, 1998, in: Textos
 Mat. S\'er. B, vol. 19, 1999, pp. 65-80.\vspace{-3mm}
 \bibitem{A01} C. Andr\'e, The basic character table of the unitriangular group, J. Algebra 241 (2001) 437-471.\vspace{-3mm}
 \bibitem{A02} C. Andr\'e, Basic characters of the unitriangular group (for arbitrary primes), Proc. Amer. Math. Soc. 130 (2002) 1934-1954.\vspace{-3mm}
\bibitem{DI2008} P. Diaconis, M. Isaacs, Supercharacters and superclasses for algebra groups, Trans. Amer. Math. Soc. 360 (2008) 2359-2392.\vspace{-3mm}
\bibitem{H12} Anders O. F. Hendrickson, Supercharacter theory constructions corresponding to Schur ring products, Comm. Algerba, 40 (2012) 4420-4438.\vspace{-3mm}
\bibitem{I1994} I. M. Isaacs, Character theory of finite groups, Dover, New York, 1994.\vspace{-3mm}
\bibitem{JL1993} G. James, M. Liebeck, Representations and Characters of Groups, Cambridge Univ. Press, Cambridge, UK, 1993.\vspace{-3mm}
\bibitem{L1991} T. Y. Lam, The first course in noncommutative rings, Springer-Verlag, New York, 1991.\vspace{-3mm}
 \bibitem{S1933} I. Schur,  Zur Theorie der einfach transitiven Permutationsgruppen.
Sitzungsber. Preuss. Akad. Wiss. Phys.-Math. Kl. (1933) 598-623. \vspace{-3mm}
\bibitem{Y01} N. Yan, Representation theory of the finite unipotent linear groups, Unpublished Ph.D. Thesis, Department of Mathematics, University of Pennsylvania, 2001.



\end{thebibliography}
\end{document}